\documentclass[12pt]{article}%
\usepackage{amsmath}
\usepackage{amsfonts}
\usepackage{amssymb}
\usepackage{graphicx}%
\providecommand{\U}[1]{\protect\rule{.1in}{.1in}}

\begin{document}

\title{\textbf{Superiorization and Perturbation Resilience of Algorithms: A
Continuously Updated Bibliography}}
\author{Yair Censor\\Department of Mathematics\\University of Haifa\\Mt.\ Carmel, Haifa 3498838, Israel\\(yair@math.haifa.ac.il)}
\date{Original report: June 13, 2015 contained 41 items.\\
First revision: March 9, 2017 contained 64 items.\\
Second revision: March 8, 2018 contained 76 items.\\
Third revision: March 11, 2019 contained 90 items.\\
Fourth revision: March 16, 2020 contained 112 items.\\
Fifth revision: March 18, 2021 contained 132 items.\\
Sixth revision: March 13, 2022 contains 150 items.\\
Seventh revision: March 14, 2023 contains 166 items.\\
Eighth revision: March 18, 2024 contains 178 items.\\
Ninth revision: March 22, 2025 contains 193 items.\\
Tenth revision: March 16, 2026 contains 205 items.}
\maketitle

\begin{abstract}
This document presents a (mostly) chronologically-ordered bibliography of
scientific publications on the superiorization methodology and perturbation
resilience of algorithms which is compiled and continuously updated by us at: http://math.haifa.ac.il/yair/bib-superiorization-censor.html.

Since the beginnings of this topic we try to trace the work that has been
published about it since its inception. To the best of our knowledge this
bibliography represents all available publications on this topic to date, and
while the URL is continuously updated we will revise this document and bring
it up to date on arXiv approximately once a year. Abstracts of the cited
works, and some links and downloadable files of preprints or reprints are
available on the above mentioned Internet page. If you know of a related
scientific work in any form that should be included here kindly write to me
on: yair@math.haifa.ac.il with full bibliographic details, a DOI if available,
and a PDF copy of the work if possible. The Internet page was initiated on
March 7, 2015, and has been last updated on March 15, 2026.

\textbf{Comment}: Some of the items have on the above mentioned Internet page
more information and links than in this report.

\textbf{Acknowledgments}: The compilation of this report was supported by the
ISF-NSFC joint research program grant No. 2874/19, by the U.S. National
Institutes of Health grant No. R01CA266467 and by the Cooperation Program in
Cancer Research of the German Cancer Research Center (DKFZ) and Israel's
Ministry of Innovation, Science and Technology (MOST).\newline\newline

\end{abstract}

\section{Trailer}

We replace the text that appeared in this section in the first version of the
report with a quotation from the preface to the special issue: Y. Censor, G.T.
Herman and M. Jiang (Guest Editors), \textquotedblleft\textit{Superiorization:
Theory and Applications}\textquotedblright, Special Issue of the journal
\textit{Inverse Problems}, Volume \textbf{33}, Number 4, April 2017
[50]\footnote{All references refer to the bibliography in the next section of
this report.}, followed by some additional notes.

\begin{quotation}
\textquotedblleft The superiorization methodology is used for improving the
efficacy of iterative algorithms whose convergence is resilient to certain
kinds of perturbations. Such perturbations are designed to `force' the
perturbed algorithm to produce more useful results for the intended
application than the ones that are produced by the original iterative
algorithm. The perturbed algorithm is called the `superiorized version' of the
original unperturbed algorithm. If the original algorithm is computationally
efficient and useful in terms of the application at hand and if the
perturbations are simple and not expensive to calculate, then the advantage of
this method is that, for essentially the computational cost of the original
algorithm, we are able to get something more desirable by steering its
iterates according to the designed perturbations.

This is a very general principle that has been used successfully in some
important practical applications, especially for inverse problems such as
image reconstruction from projections, intensity-modulated radiation therapy
and nondestructive testing, and awaits to be implemented and tested in
additional fields.

An important case is when the original algorithm is `feasibility-seeking' (in
the sense that it strives to find some point that is compatible with a family
of constraints) and the perturbations that are introduced into the original
iterative algorithm aim at reducing (not necessarily minimizing) a given merit
function. In this case superiorization has a unique place in optimization
theory and practice.

Many constrained optimization methods are based on methods for unconstrained
optimization that are adapted to deal with constraints. Such is, for example,
the class of projected gradient methods wherein the unconstrained minimization
inner step `leads' the process and a projection onto the whole constraints set
(the feasible set) is performed after each minimization step in order to
regain feasibility. This projection onto the constraints set is in itself a
non-trivial optimization problem and the need to solve it in every iteration
hinders projected gradient methods and limits their efficiency to only
feasible sets that are `simple to project onto.' Barrier or penalty methods
likewise are based on unconstrained optimization combined with various
`add-on's that guarantee that the constraints are preserved. Regularization
methods embed the constraints into a `regularized' objective function and
proceed with unconstrained solution methods for the new regularized objective function.

In contrast to these approaches, the superiorization methodology can be viewed
as an antipodal way of thinking. Instead of adapting unconstrained
minimization algorithms to handling constraints, it adapts feasibility-seeking
algorithms to reduce merit function values. This is done while retaining the
feasibility-seeking nature of the algorithm and without paying a high
computational price. Furthermore, general-purpose approaches have been
developed for automatically superiorizing iterative algorithms for large
classes of constraints sets and merit functions; these provide algorithms for
many application tasks.\textquotedblright
\end{quotation}

To a novice on the superiorization methodology and perturbation resilience of
algorithms we recommend to read first the recent reviews in [16, 25, 39, 162].
For a recent description of previous work that is related to superiorization
but is not included here, such as the works of Sidky and Pan, e.g., [6], we
direct the reader to [24, section~3]. The SNARK14 software package [42], with
its in-built capability to superiorize iterative algorithms to improve their
performance, can be helpful to practitioners. Naturally there is variability
among the bibliography items below in their degree of relevance to the
superiorization methodology and perturbation resilience of algorithms. In
some, such as in, e.g., [23] below, superiorization does not appear in the
title, abstract or introduction but only inside the work, e.g., [23,
Subsection 6.2.1: Optimization vs. Superiorization].

A word about the history. The terms and notions \textquotedblleft
superiorization\textquotedblright\ and \textquotedblleft perturbation
resilience\textquotedblright\ first appeared in the 2009 paper of Davidi,
Herman and Censor [7] which followed its 2007 forerunner by Butnariu, Davidi,
Herman and Kazantsev [3]. The ideas have some of their roots in the 2006 and
2008 papers of Butnariu, Reich and Zaslavski [2, 4]. All these culminated in
Ran Davidi's 2010 Ph.D. dissertation [13].

\section{The Bibliography}

[1] P.L. Combettes, On the numerical robustness of the parallel projection
method in signal synthesis, \textit{IEEE Signal Processing Letters}, Vol. 8,
pp. 45--47, (2001). DOI:10.1109/97.895371.\medskip

[2] D. Butnariu, S. Reich and A.J. Zaslavski, Convergence to fixed points of
inexact orbits of Bregman-monotone and of nonexpansive operators in Banach
spaces, in: H.F. Nathansky, B.G. de Buen, K. Goebel, W.A. Kirk, and B. Sims
(Editors), \textit{Fixed Point Theory and its Applications}, (Conference
Proceedings, Guanajuato, Mexico, 2005), Yokahama Publishers, Yokahama, Japan,
pp. 11--32, 2006.
http://www.ybook.co.jp/pub/ISBN\%20978-4-9465525-0.htm.\medskip

[3] D. Butnariu, R. Davidi, G.T. Herman and I.G. Kazantsev, Stable convergence
behavior under summable perturbations of a class of projection methods for
convex feasibility and optimization problems, \textit{IEEE Journal of Selected
Topics in Signal Processing}, Vol. 1, pp. 540--547, (2007).
DOI:10.1109/JSTSP.2007.910263.\medskip

[4] D. Butnariu, S. Reich and A.J. Zaslavski, Stable convergence theorems for
infinite products and powers of nonexpansive mappings, \textit{Numerical
Functional Analysis and Optimization}, Vol. 29, pp. 304--323, (2008).
DOI:10.1080/01630560801998161.\medskip

[5] G.T. Herman and R. Davidi, Image reconstruction from a small number of
projections, \textit{Inverse Problems}, Vol. 24, 045011 (17pp),
(2008).\newline DOI:10.1088/0266-5611/24/4/045011.\medskip

[6] E.Y. Sidky and X. Pan, Image reconstruction in circular cone-beam computed
tomography by constrained, total-variation minimization, \textit{Physics in
Medicine and Biology}, Vol. 53, pp. 4777--4807, (2008).
DOI:10.1088/0031-9155/53/17/021.\medskip

[7] R. Davidi, G.T. Herman and Y. Censor, Perturbation-resilient
block-iterative projection methods with application to image reconstruction
from projections, \textit{International Transactions in Operational Research},
Vol. 16, pp. 505--524, (2009). DOI:10.1111/j.1475-3995.2009.00695.x.\medskip

[8] G.T. Herman, \textit{Fundamentals of Computerized Tomography: Image
Reconstruction from Projections}, Springer-Verlag, London, UK, 2nd Edition,
2009. DOI:10.1007/978-1-84628-723-7.\medskip

[9] S.N. Penfold, \textit{Image Reconstruction and Monte Carlo Simulations in
the Development of Proton Computed Tomography for Applications in Proton
Radiation Therapy}, Ph.D. dissertation, Centre for Medical Radiation Physics,
University of Wollongong, 2010.\newline
http://ro.uow.edu.au/cgi/viewcontent.cgi?article=4305\&context=theses.\medskip

[10] S.N. Penfold, R.W. Schulte, Y. Censor, V. Bashkirov, S. McAllister, K.E.
Schubert and A.B. Rosenfeld, Block-iterative and string-averaging projection
algorithms in proton computed tomography image reconstruction, in: Y. Censor,
M. Jiang and G. Wang (Editors), \textit{Biomedical Mathematics: Promising
Directions in Imaging, Therapy Planning and Inverse Problems}, Medical Physics
Publishing, Madison, WI, USA, 2010, pp. 347--367.
https://www.medicalphysics.org/SimpleCMS.php?content=reviewlist.php\&\newline
isbn=9781930524484.\medskip

[11] Y. Censor, R. Davidi and G.T. Herman, Perturbation resilience and
superiorization of iterative algorithms, \textit{Inverse Problems}, Vol. 26,
(2010) 065008 (12pp). DOI:10.1088/0266-5611/26/6/065008.\medskip

[12] S.N. Penfold, R.W. Schulte, Y. Censor and A.B. Rosenfeld, Total variation
superiorization schemes in proton computed tomography image reconstruction,
\textit{Medical Physics}, Vol. 37, pp. 5887--5895, (2010).\newline
DOI:10.1118/1.3504603.\medskip

[13] R. Davidi, \textit{Algorithms for Superiorization and their Applications
to Image Reconstruction}, Ph.D. dissertation, Department of Computer Science,
The City University of New York, NY, USA, 2010.\newline
http://gradworks.umi.com/34/26/3426727.html.\medskip

[14] E. Gardu\~{n}o, R. Davidi and G.T. Herman, Reconstruction from a few
projections by $\ell_{1}$-minimization of the Haar transform, \textit{Inverse
Problems}, Vol. 27, 055006, (2011). DOI:10.1088/0266-5611/27/5/055006.\medskip

[15] Y. Censor, W. Chen, P.L. Combettes, R. Davidi and G.T. Herman, On the
effectiveness of projection methods for convex feasibility problems with
linear inequality constraints, \textit{Computational Optimization and
Applications}, Vol. 51, pp. 1065--1088, (2012). DOI:10.1007/s10589-011-9401-7.
A related (unpublished) Technical Report: W. Chen, Data sets of very large
linear feasibility problems solved by projection methods, March 2,
2011.\medskip

[16] G.T. Herman, E. Gardu\~{n}o, R. Davidi and Y. Censor, Superiorization: An
optimization heuristic for medical physics, \textit{Medical Physics}, Vol. 39,
pp. 5532--5546, (2012). DOI:10.1118/1.4745566.\medskip

[17] R. Davidi, R.W. Schulte, Y. Censor and L. Xing, Fast superiorization
using a dual perturbation scheme for proton computed tomography,
\textit{Transactions of the American Nuclear Society}, Vol. 106, pp. 73--76,
(2012).\medskip

[18] T. Nikazad, R. Davidi and G.T. Herman, Accelerated perturbation-resilient
block-iterative projection methods with application to image reconstruction,
\textit{Inverse Problems}, Vol. 28, 035005 (19pp), (2012).\newline
DOI:10.1088/0266-5611/28/3/035005.\medskip

[19] D. Steinberg, V. Bashkirov, V. Feng, R.F. Hurley, R.P. Johnson, S.
Macafee, T. Plautz, H.F.-W. Sadrozinski, R. Schulte and A. Zatserklyaniy,
Monte Carlo simulations for the development a clinical proton CT scanner,
\textit{Nuclear Science Symposium and Medical Imaging Conference (NSS/MIC)},
2012 IEEE, pp. 1311--1315. Oct. 27-Nov. 3, 2012, Anaheim, CA, USA.\newline
DOI:10.1109/NSSMIC.2012.6551320.\medskip

[20] W. Jin, Y. Censor and M. Jiang, A heuristic superiorization-like approach
to bioluminescence, \textit{International Federation for Medical and
Biological Engineering (IFMBE) Proceedings}, Vol. 39, pp. 1026--1029, (2013).
DOI:10.1007/978-3-642-29305-4\_269.\medskip

[21] Y. Censor and A.J. Zaslavski, Convergence and perturbation resilience of
dynamic string-averaging projection methods, \textit{Computational
Optimization and Applications}, Vol. 54, pp. 65--76, (2013).
DOI:10.1007/s10589-012-9491-x.\medskip

[22] S.-S. Luo, \textit{Reconstruction Algorithms for Single-photon Emission
Computed Tomography}, Ph.D. dissertation, Computational Mathematics, \newline
Peking University (PKU), Beijing, P.R. China, 2013.\newline
http://www.dissertationtopic.net/doc/2220625.\medskip

[23] X. Zhang, \textit{Prior-Knowledge-Based Optimization Approaches for CT
Metal Artifact Reduction}, Ph.D. dissertation, Dept. of Electrical
Engineering, Stanford University, Stanford, CA, USA, 2013.\newline
http://purl.stanford.edu/ws303zb5770.\medskip

[24] Y. Censor, R. Davidi, G.T. Herman, R.W. Schulte and L. Tetruashvili,
Projected subgradient minimization versus superiorization, \textit{Journal of
Optimization Theory and Applications}, Vol. 160, pp. 730--747, (2014).
DOI:10.1007/s10957-013-0408-3.\medskip

[25] G.T. Herman, Superiorization for image analysis, in:
\textit{Combinatorial Image Analysis}, Lecture Notes in Computer Science Vol.
8466, Springer, 2014, pp. 1--7. DOI:10.1007/978-3-319-07148-0\_1.\medskip

[26] S. Luo and T. Zhou, Superiorization of EM algorithm and its application
in single-photon emission computed tomography (SPECT), \textit{Inverse
Problems and Imaging}, Vol. 8, pp. 223--246, (2014). \newline
DOI:10.3934/ipi.2014.8.223.\medskip

[27] M.J. Schrapp and G.T. Herman, Data fusion in X-ray computed tomography
using a superiorization approach, \textit{Review of Scientific Instruments},
Vol. 85, 053701 (9pp), (2014). DOI:10.1063/1.4872378.\medskip

[28] M. Schrapp, M. Goldammer, K. Sch\"{o}rner and J. Stephan, Improvement of
image quality in computed tomography via data fusion, \textit{Proceedings of
the 5th International Conference on Industrial Computed Tomography (iCT)}, pp.
283--289, February 2014, the University of Applied Sciences, Wels, Upper
Austria. http://www.ndt.net/article/ctc2014/papers/283.pdf.\medskip

[29] E. Gardu\~{n}o and G.T. Herman, Superiorization of the ML-EM algorithm,
\textit{IEEE Transactions on Nuclear Science}, Vol. 61, pp. 162--172, (2014).
DOI:10.1109/TNS.2013.2283529.\medskip

[30] O. Langthaler, \textit{Incorporation of the Superiorization Methodology
into Biomedical Imaging Software}, Marshall Plan Scholarship Report, Salzburg
University of Applied Sciences, Salzburg, Austria, and The Graduate Center of
the City University of New York, NY, USA, September 2014, (76 pages).
http://www.marshallplan.at/images/papers\_ scholarship/2014/Salzburg\_\newline
University\_of\_Applied\_Sciences\_LangthalerOliver\_2014.pdf.\medskip

[31] B. Prommegger, \textit{Verification and Evaluation of Superiorized
Algorithms Used in Biomedical Imaging: Comparison of Iterative Algorithms With
and Without Superiorization for Image Reconstruction from Projections},
Marshall Plan Scholarship Report, Salzburg University of Applied Sciences,
Salzburg, Austria, and The Graduate Center of the City University of New York,
NY, USA, October 2014, (84 pages).\newline
http://www.marshallplan.at/images/papers\_scholarship/2014/Salzburg\_\newline
University\_of\_Applied\_Sciences\_PrommeggerBernhard\_2014.pdf.\medskip

[32] D.C. Hansen, \textit{Improving Ion Computed Tomography}, Ph.D.
dissertation, Aarhus University, Experimental Clinical Oncology, Aarhus,
Denmark, 2014.
http://pure.au.dk//portal/files/83515131/dissertation.pdf.\medskip

[33] J. Lee, C. Kim, B. Min, J. Kwak, S. Park, S-B. Lee, S. Park and S. Cho,
Sparse-view proton computed tomography using modulated proton beams,
\textit{Medical Physics}, Vol. 42, pp. 1129--1137, (2015).\newline
DOI:10.1118/1.4906133.\medskip

[34] T. Nikazad and M. Abbasi, Perturbation-resilient iterative methods with
an infinite pool of mappings, \textit{SIAM Journal on Numerical Analysis},
Vol. 53, pp. 390--404, (2015). DOI:10.1137/14095724X.\medskip

[35] F. Arroyo, E. Arroyo, X. Li and J. Zhu, The convergence of the block
cyclic projection with an overrelaxation parameter for compressed sensing
based tomography, \textit{Journal of Computational and Applied Mathematics},
Vol. 280, pp. 59--67, (2015). DOI:10.1016/j.cam.2014.11.036.\medskip

[36] R. Davidi, Y. Censor, R.W. Schulte, S. Geneser and L. Xing,
Feasibility-seeking and superiorization algorithms applied to inverse
treatment planning in radiation therapy, \textit{Contemporary Mathematics},
Vol. 636, pp. 83--92, (2015). DOI:10.1090/conm/636/12729.\medskip

[37] Y. Censor and D. Reem, Zero-convex functions, perturbation resilience,
and subgradient projections for feasibility-seeking methods,
\textit{Mathematical Programming, Series A}, Vol. 152, pp. 339--380,
(2015).\newline DOI:10.1007/s10107-014-0788-7.\medskip

[38] Y. Censor and A.J. Zaslavski, Strict Fej\'{e}r monotonicity by
superiorization of feasibility-seeking projection methods, \textit{Journal of
Optimization Theory and Applications}, Vol. 165, pp. 172--187, (2015).
DOI:10.1007/s10957-014-0591-x.\medskip

[39] Y. Censor, Weak and strong superiorization: Between feasibility-seeking
and minimization, \textit{Analele Stiintifice ale Universitatii Ovidius
Constanta, Seria Matematica}, Vol. 23, pp. 41--54, (2015).
DOI:10.1515/auom-2015-0046.\medskip

[40] H.H. Bauschke and V.R. Koch, Projection methods: Swiss army knives for
solving feasibility and best approximation problems with half-spaces,
\textit{Contemporary Mathematics}, Vol. 636, pp. 1--40, (2015).\newline
DOI:10.1090/conm/636/12726.\medskip

[41] M.J. Schrapp, \textit{Multi Modal Data Fusion in Industrial X-ray
Computed Tomography}, Ph.D. dissertation, Fakult\"{a}t f\"{u}r Physik der
Technischen Universit\"{a}t M\"{u}nchen, Munich, Germany, 2015.\medskip

[42] SNARK14, \textit{A programming system for the reconstruction of 2D images
from 1D projections designed to help researchers in developing and evaluating
reconstruction algorithms}. In particular, SNARK14 can be used for automatic
superiorization of any iterative reconstruction algorithm. Released:
2015.\medskip

[43] W. Jin, Y. Censor and M. Jiang, Bounded perturbation resilience of
projected scaled gradient methods, C\textit{omputational Optimization and
Applications}, Vol. 63, pp. 365--392, (2016).
DOI:10.1007/s10589-015-9777-x.\medskip

[44] Q-L. Dong, J. Zhao and S. He, Bounded perturbation resilience of the
viscosity algorithm, \textit{Journal of Inequalities and Applications},
2016:299 (12pp), 2016. DOI:10.1186/s13660-016-1242-6.\medskip

[45] E. Nurminski, Finite-value superiorization for variational inequality
problems, arXiv:1611.09697, (2016). [arXiv:1611.09697].\medskip

[46] S. Luo, Y. Zhang, T. Zhou and J. Song, Superiorized iteration based on
proximal point method and its application to XCT image reconstruction,
arXiv:1608.03931, (2016). [arXiv:1608.03931].\medskip

[47] Y. Censor and Y. Zur, Linear superiorization for infeasible linear
programming, in: Y. Kochetov, M. Khachay, V. Beresnev, E. Nurminski and P.
Pardalos (Editors), \textit{Discrete Optimization and Operations Research},
Lecture Notes in Computer Science (LNCS), Vol. 9869, (2016), Springer
International Publishing, pp. 15--24. DOI:10.1007/978-3-319-44914-2\_2.
Reprint of the paper is available for free download on the publisher's website
at:\newline
http://www.springer.com/gp/book/9783319449135?wt\_mc=ThirdParty.\newline
SpringerLink.3.EPR653.About\_eBook, under the link: \textquotedblleft Download
Sample pages 2 PDF (774.4 KB)\textquotedblright\ thereon.\medskip

[48] C. Havas, \textit{Revised Implementation and Empirical Study of Maximum
Likelihood Expectation Maximization Algorithms with and without
Superiorization in Image Reconstruction}, Marshall Plan Scholarship Report,
Salzburg University of Applied Sciences, Salzburg, Austria, and The Graduate
Center of the City University of New York, NY, USA, October 2016, (49
pages).\newline
https://static1.squarespace.com/static/559921a3e4b02c1d7480f8f4/t/596c97\newline
aad1758e1c6808c0fa/1500288944245/ Havas+Clemens\_615.pdf.\medskip

[49] T. Humphries, J. Winn and A. Faridani, Superiorized algorithm for
reconstruction of CT images from sparse-view and limited-angle polyenergetic
data, arXiv:1701.03396, (2017). [arXiv:1608.03931].\medskip

[50] Y. Censor, G.T. Herman and M. Jiang, Guest Editors, Superiorization:
Theory and Applications, Special Issue of the journal Inverse Problems, Volume
33, Number 4, April 2017. Read the titles and abstracts of all 14 papers
included in the special issue on the journal's website at:\newline
http://iopscience.iop.org/issue/0266-5611/33/4;jsessionid=AF091E29223E\newline%
16A29F38C99720D302B6.c4.iopscience.cld.iop.org.\medskip

[51] D. Reem and A. De Pierro, A new convergence analysis and perturbation
resilience of some accelerated proximal forward-backward algorithms with
errors, \textit{Inverse Problems}, Vol. 33 (2017), 044001.\newline
https://doi.org/10.1088/1361-6420/33/4/044001.\medskip

[52] T. Nikazad and M. Abbasi, A unified treatment of some perturbed fixed
point iterative methods with an infinite pool of operators, \textit{Inverse
Problems}, Vol. 33 (2017), 044002.\newline
https://doi.org/10.1088/1361-6420/33/4/044002.\medskip

[53] M. Yamagishi and I. Yamada, Nonexpansiveness of a linearized augmented
Lagrangian operator for hierarchical convex optimization, \textit{Inverse
Problems}, Vol. 33 (2017), 044003.\newline
https://doi.org/10.1088/1361-6420/33/4/044003.\medskip

[54] A.J. Zaslavski, Asymptotic behavior of two algorithms for solving common
fixed point problems, \textit{Inverse Problems}, Vol. 33 (2017), 044004.
https://doi.org/10.1088/1361-6420/33/4/044004.\medskip

[55] S. Reich and A.J. Zaslavski, Convergence to approximate solutions and
perturbation resilience of iterative algorithms, \textit{Inverse Problems},
Vol. 33 (2017), 044005. https://doi.org/10.1088/1361-6420/33/4/044005.\medskip

[56] Y. Censor, Can linear superiorization be useful for linear optimization
problems?, \textit{Inverse Problems}, Vol. 33 (2017), 044006.\newline
https://doi.org/10.1088/1361-6420/33/4/044006.\medskip

[57] H. He and H-K. Xu, Perturbation resilience and superiorization
methodology of averaged mappings, \textit{Inverse Problems}, Vol. 33 (2017),
044007. \newline https://doi.org/10.1088/1361-6420/33/4/044007.\medskip

[58] H-K. Xu, Bounded perturbation resilience and superiorization techniques
for the projected scaled gradient method, \textit{Inverse Problems}, Vol. 33
(2017), 044008. https://doi.org/10.1088/1361-6420/33/4/044008.\medskip

[59] A. Cegielski and F. Al-Musallam, Superiorization with level control,
\textit{Inverse Problems}, Vol. 33 (2017), 044009.\newline
https://doi.org/10.1088/1361-6420/aa5d79.\medskip

[60] E.S. Helou, M.V.W. Zibetti and E.X. Miqueles, Superiorization of
incremental optimization algorithms for statistical tomographic image
reconstruction, \textit{Inverse Problems}, Vol. 33 (2017), 044010.\newline
https://doi.org/10.1088/1361-6420/33/4/044010.\medskip

[61] E. Gardu\~{n}o and G.T. Herman, Computerized tomography with total
variation and with shearlets, \textit{Inverse Problems}, Vol. 33 (2017),
044011. https://doi.org/10.1088/1361-6420/33/4/044011.\medskip

[62] E. Bonacker, A. Gibali, K-H. K\"{u}fer and P. S\"{u}ss, Speedup of
lexicographic optimization by superiorization and its applications to cancer
radiotherapy treatment, \textit{Inverse Problems}, Vol. 33 (2017),
044012.\newline https://doi.org/10.1088/1361-6420/33/4/044012.\medskip

[63] J. Zhu and S. Penfold, Total variation superiorization in dual-energy CT
reconstruction for proton therapy treatment planning, \textit{Inverse
Problems}, Vol. 33 (2017), 044013.
https://doi.org/10.1088/1361-6420/33/4/044013.\medskip

[64] Q. Yang, W. Cong and G. Wang, Superiorization-based multi-energy CT image
reconstruction, \textit{Inverse Problems}, Vol. 33 (2017), 044014.
https://\newline doi.org/10.1088/1361-6420/aa5e0a.\medskip

[65] T. Nikazad, M. Abbasi and T. Elfving, Error minimizing relaxation
strategies in Landweber and Kaczmarz type iterations, \textit{Journal of
Inverse and Ill-posed Problems}, Vol. 25, pp. 35--56, (2017).\newline
DOI:10.1515/jiip-2015-0082.\medskip

[66] Q.L. Dong, Y.J. Cho and Th.M. Rassias, Multi-step inertial
Krasno-\newline sel'ski-Mann algorithm for nonexpansive operators,
\textit{Technical Report. }Preprint from ResearchGate at:
https://www.researchgate.net/publication/318440732\newline%
\_MiKM\_Multi-step\_inertial\_Krasnosel\%27skiiskii-Mann\_algorithm\_and\_\newline
its\_applications. (2017).\medskip

[67] C.O.S. Sorzano, J. Vargas, J. Ot\'{o}n, J.M. de la Rosa-Trev\'{\i}n, J.L.
Vilas, M. Kazemi, R. Melero, L. del Ca\~{n}o, J. Cuenca, P. Conesa, J.
G\'{o}mez-Blanco, R. Marabini and J.M. Carazo, A survey of the use of
iterative reconstruction algorithms in electron microscopy, \textit{BioMed
Research International}, Vol. 2017 (2017), Article ID 6482567, 17 pages,
https://doi.org/10.1155/2017\newline/6482567.\medskip

[68] C.L. Byrne, The Dykstra and Bregman-Dykstra algorithms as superiorization
(December 26, 2017). \textit{Technical report}, (2017). Preprint available
from ResearchGate.\medskip

[69] Q.-L. Dong, A. Gibali, D. Jiang and S.-H. Ke, Convergence of projection
and contraction algorithms with outer perturbations and their applications to
sparse signals recovery, \textit{Fixed Point Theory and Applications}, (2018)
20: 16. https://doi.org/10.1007/s11784-018-0501-1. \medskip

[70] Q.-L. Dong, A. Gibali, D. Jiang and Y.-C. Tang, Bounded perturbation
resilience of extragradient-type methods and their applications,
\textit{Journal of Inequalities and Applications}, (2017) 2017:280.
DOI:10.1186/s13660-017-1555-0.\medskip

[71] A. Gibali and S. Petra, DC-programming versus $\ell_{0}$-superiorization
for discrete tomography, \textit{Analele Stiintifice ale Universitatii Ovidius
Constanta-Seria Matematica,} Vol. 26, 2018, pp. 105-133. Available at the
journal's homepage at: \newline
http://www.anstuocmath.ro/mathematics//Anale2018vol2/6\_GibaliA.pdf.\medskip

[72] Yanni Guo, W. Cui and Yansha Guo, Perturbation resilience of proximal
gradient algorithm for composite objectives, \textit{Journal of Nonlinear
Sciences and Applications (JNSA)}, Vol. 10, pp. 5566-5575, (2017),
http://dx.doi.org/10.22436/jnsa.010.10.36.\medskip

[73] C. Bargetz, S. Reich and R. Zalas, Convergence properties of dynamic
string averaging projection methods in the presence of perturbations,
\textit{Numerical Algorithms}, Vol. 77, pp. 185-209, (2018).
https://doi.org/10.1007/s11075-017-0310-4.\medskip

[74] M.V.W. Zibetti, C. Lin and G.T. Herman, Total variation superiorized
conjugate gradient method for image reconstruction, \textit{Inverse Problems},
Vol. 34 (2018), 034001. https://doi.org/10.1088/1361-6420/aaa49b.\medskip

[75] A. Gibali, K-H. K\"{u}fer, D. Reem and P. S\"{u}ss, A generalized
projection-based scheme for solving convex constrained optimization problems,
\textit{Computational Optimization and Applications}, Vol. 70, pp. 737-762,
(2018). https://doi.org/10.1007/s10589-018-9991-4. \medskip

[76] B. Schultze, Y. Censor, P. Karbasi, K.E. Schubert, and R.W. Schulte, An
improved method of total variation superiorization applied to reconstruction
in proton computed tomography, \textit{IEEE Transactions on Medical Imaging},
accepted for publication. (2019). Available on arXiv at: \newline
https://arxiv.org/abs/1803.01112.\medskip

[77] Y. Censor, H. Heaton, and R.W. Schulte, Derivative-free superiorization
with component-wise perturbations. \textit{Numerical Algorithms}, accepted for
publication. (2018). https://doi.org/10.1007/s11075-018-0524-0.\medskip

[78] Y. Guo and W. Cui, Strong convergence and bounded perturbation resilience
of a modified proximal gradient algorithm, \textit{Journal of Inequalities and
Applications}, 2018:103, (2018), https://doi.org/10.1186/s13660-018-1695-x.
Free download of full paper from the publisher at:\newline
https://journalofinequalitiesandapplications.springeropen.com/track/pdf/\newline%
10.1186/s13660-018-1695-x.\medskip

[79] A.J. Zaslavski, \textit{Algorithms for Solving Common Fixed Point
Problems}, Springer International Publishing AG, part of Springer Nature,
(2018). Freely available from the publisher at: \newline
https://link.springer.com/book/ 10.1007\%2F978-3-319-77437-4.\medskip

[80] T.Y. Kong, H. Pajoohesh and G.T. Herman, String-averaging algorithms for
convex feasibility with infinitely many sets, June (2018). Preprint from arXiv
at: https://arxiv.org/abs/1807.00234.\medskip

[81] T. Nikazad and M. Abbasi, Perturbed fixed point iterative methods based
on pattern structure, \textit{Mathematical Methods in the Applied Sciences},
(2018), 1-11. https://doi.org/10.1002/mma.5100. Free download of full paper
from the publisher at: \newline
https://onlinelibrary.wiley.com/doi/epdf/10.1002/mma.5100.\medskip

[82] E.S. Helou, G.T. Herman, C. Lin and M.V.W. Zibetti, Superiorization of
preconditioned conjugate gradient algorithms for tomographic image
reconstruction, \textit{Applied Analysis and Optimization}, Vol. 2, pp.
271-284, (2018). Readable on Yokohama Publishers webpage at: \newline
http://www.ybook.co.jp/online2/aaov2n2.html.\medskip

[83] S. Luo, Y. Zhang, T. Zhou, J. Song and Y. Wang, XCT image reconstruction
by a modified superiorized iteration and theoretical analysis,
\textit{Optimization Methods and Software}, Published online: 02 Jan 2019.
(2019). Readable at https://doi.org/10.1080/10556788.2018.1560442.\medskip

[84] P. Duan, X. Zheng and J. Zhao, Strong convergence theorems of viscosity
iterative algorithms for split common fixed point problems,
\textit{Mathematics}, 2019, 7(1), 14;(2019). Open access at:
https://www.mdpi.com/2227-7390/7/1/14/htm.\medskip

[85] Q.L. Dong, J. Huang, X.H. Li, Y.J. Cho and Th.M. Rassias, MiKM:
multi-step inertial Krasnosel'ski\u{\i}-Mann algorithm and its applications,
\textit{Journal of Global Optimization}, accepted for publication,
(2019).\newline
https://link.springer.com/article/10.1007\%2Fs10898-018-0727-x.\medskip

[86] R. Davidi, R.M. Haralick and G.T. Herman, Derivative-free superiorization
using the facet model, a presentation at the Seminar on Image Processing and
Computer Vision, September 16, 2010. The Graduate Center, City University of
New York (CUNY) (2010). Get the presentation slides at:
http://math.haifa.ac.il/yair/WEBPAGE-AFTER-200915/Presentation\_\newline
Seminar\_8.pdf. A note from the page maintainer: In spite of this document
being from 2010, I place it here because it came to my attention only very
recently. y.c.\medskip

[87] M.A. Kalkhoran and D. Vray, Sparse sampling and reconstruction for an
optoacoustic ultrasound volumetric hand-held probe, \textit{Biomedical Optics
Express}, Vol. 10, pp. 1545-1556, (2019). Accessible at OSA Publishing webpage
at: https://www.osapublishing.org/boe/abstract.cfm?uri=boe-10-4-1545.\medskip

[88] T. Humphries, PSARTSUP: a GitHub archive that contains code used for the
Superiorized pSART method described in: \textquotedblleft Superiorized
algorithm for reconstruction of CT images from sparse-view and limited-angle
polyenergetic data\textquotedblright\ by Humphries, Winn and Faridani here,
see item [49] above, and in: \textquotedblleft Superiorized polyenergetic
reconstruction algorithm for reduction of metal artifacts in CT
images\textquotedblright\ by Humphries and Gibali here, see item [89] below.
Link to the PSARTSUP GitHub archive at: \newline
https://github.com/TDHumphries/PSARTSUP. (2017). A note from the page
maintainer: In spite of this document being from 2017, I place it here because
it came to my attention only very recently. y.c.\medskip

[89] A. Gibali and T. Humphries, Superiorized polyenergetic reconstruction
algorithm for reduction of metal artifacts in CT images. \textit{Proceedings
of the 2017 IEEE Nuclear Science Symposium and Medical Imaging Conference
(NSS/MIC 2017)}, October 21-28, 2017, Atlanta, Georgia, USA, pp. 920-925.
(2017). See Proceedings table of contents at: \newline
http://toc.proceedings.com/41689webtoc.pdf. A note from the page maintainer:
In spite of this document being from 2017, I place it here because it came to
my attention only very recently. y.c.\medskip

[90] G.T. Herman, Iterative reconstruction techniques and their
superiorization for the inversion of the Radon transform, in: R. Ramlau and O.
Scherzer, eds., \textit{The Radon Transform: The First 100 Years and Beyond},
Chapter 10, pp. 217-238 (2019), De Gruyter, Berlin, Boston.
https://doi.org/\newline10.1515/9783110560855-010.\medskip

[91] R. Cassetta, P. Piersimoni, M. Riboldi, V. Giacometti, V. Bashkirov, G.
Baroni, C. Ordonez, G. Coutrakon and R. Schulte, Accuracy of low-dose proton
CT image registration for pretreatment alignment verification in reference to
planning proton CT, \textit{Journal of Applied Clinical Medical Physics},
(2019);20:4:83-90. https://doi.org/10.1002/acm2.12565.\medskip

[92] J. Fink, R.L.G. Cavalcante and S. Stanczak, Multicast Beamforming Using
Semidefinite Relaxation and Bounded Perturbation Resilience, \textit{2019 IEEE
International Conference on Acoustics, Speech and Signal Processing} (ICASSP),
(2019). DOI:10.1109/ICASSP.2019.8682325.\medskip

[93] M. Hoseini, Superiorization and its importance in the optimization,
\textit{ICNS Conference Proceeding}, Vol. 3. Issue 2, pp. 6-9. (2019).
Available at Charmo University, Kurdistan Region-Iraq.
http://dx.doi.org/10.31530/17030.\medskip

[94] E. Bonacker, A. Gibali and K.-H. K\"{u}fer, Accelerating Two Projection
Methods via Perturbations with Application to Intensity-Modulated Radiation
Therapy, \textit{Applied Mathematics \& Optimization}, (2019), accepted for
publication. https://doi.org/10.1007/s00245-019-09571-4.\medskip

[95] C.L. Byrne, Thoughts on superiorization. This is a sequence of five short
preprints posted on ResearchGate between April 28 and June 21, 2019. A
combined PDF is Here.\medskip

[96] J. Fink, D. Schaeufele, M. Kasparick, R.L.G. Cavalcante and S. Stanczak,
Cooperative Localization by Set-theoretic Estimation, WSA 2019; \textit{23rd
International ITG Workshop on Smart Antennas}, pp. 1-8. WSA 2019, April 24-26,
2019, Vienna, Austria. (2019). Open access at: IEEE Xplore Digital
Library.\medskip

[97] S. Reich and R. Zalas, A modular string averaging procedure for solving
the common fixed point problem for quasi-nonexpansive mappings in Hilbert
space, \textit{Numerical Algorithms}, Vol. 72, pp. 297-323, (2016).
https://doi.org/10.1007/s11075-015-0045-z. Reprint readable on: Springer
Nature SharedIt. A note from the page maintainer: Bounded perturbation
resilience appears in Example 4.7. y.c.\medskip

[98] Yanni Guo and Xiaozhi Zhao, Bounded Perturbation Resilience and
Superiorization of Proximal Scaled Gradient Algorithm with Multi-Parameters,
\textit{Mathematics}, 7(6), 535 (14 pp) (2019).
https://doi.org/10.3390/math7060535. Open access at: Scholarly open access
publishing, MDPI.\medskip

[99] S. Rai, Image Quality Measures in Proton Computed Tomography, M.Sc.
(Master of Science) Thesis, Department of Physics, Northern Illinois
University (NIU), De Kalb, IL, USA. (110 pp) (2015).\newline
https://commons.lib.niu.edu/handle/10843/18827. Available at: Northern
Illinois University's digital repository .\medskip

[100] Y. Censor, E. Gardu\~{n}o, E.S. Helou and G.T. Herman, Derivative-Free
Superiorization: Principle and Algorithm, \textit{Preprint} (2019). \newline
https://arxiv.org/abs/1908.10100.\medskip

[101] Y. Censor, E. Levy, An analysis of the superiorization method via the
principle of concentration of measure, \textit{Applied Mathematics and
Optimization}, accepted for publication (2019).
https://arxiv.org/abs/1909.00398.\medskip

[102] C.L. Byrne, What do simulations tell us about superiorization?
\textit{Preprint posted on ResearchGate}. September 2019.\medskip

[103] E. Bonacker, \textit{Perturbed Projection Methods in Convex Optimization
- Applied to Radiotherapy Planning}, Ph.D. Dissertation, Fachbereich
Mathematik der Technischen Universit\"{a}t Kaiserslautern, Kaiserslautern,
Germany. July 2019.\medskip

[104] Peichao Duan and Xubang Zheng, Bounded perturbation resilience and
superiorization techniques for a modified proximal gradient method,
\textit{Optimization}. (2019). [Reprint from
https://www.tandfonline.com/doi/full/\newline%
10.1080/02331934.2019.1686631].\medskip

[105] Mohsen Hoseini, Shahram Saeidi and Do Sang Kim, On perturbed hybrid
steepest descent method with minimization or superiorization for
subdifferentiable functions, \textit{Numerical Algorithms}, Vol. 85, pp.
353-374, (2020). https://doi.org/10.1007/s11075-019-00818-3.\medskip

[106] Yamin Wang, Fenghui Wang and Haixia Zhang, Strong convergence of
viscosity forward-backward algorithm to the sum of two accretive operators in
Banach space, \textit{Optimization}, (2019).
https://doi.org/10.1080/02331934.2019.\newline1705299. Full access at:
https://www.tandfonline.com/.\medskip

[107] Yanni Guo and Xiaozhi Zhao, Strong convergence of over-relaxed
multiparameter proximal scaled gradient algorithm and superiorization,
\textit{Optimization}, (2020). https://doi.org/10.1080/02331934.2020.1722124.
Full access at: https://www.tandfonline.com/.\medskip

[108] Nuttapol Pakkaranang, Poom Kumam, Vasile Berinde and Yusuf I. Suleiman,
Superiorization methodology and perturbation resilience of inertial proximal
gradient algorithm with application to signal recovery, \textit{The Journal of
Supercomputing}, (2020). https://doi.org/10.1007/s11227-020-03215-z.\medskip

[109] Thomas Humphries and Boyang (Jessie) Wang, Superiorized method for metal
artifact reduction. \textit{Medical Physics}, Vol. 47, pp. 3984-3995 (2020).
Published on line at:
https://aapm.onlinelibrary.wiley.com/doi/abs/10.1002/\newline
mp.14332.\medskip

[110] Eliahu Levy, Accumulated Random Distances in High Dimensions-Ways of
Calculation. Preprint. (2020). https://arxiv.org/abs/2003.10941.\medskip

[111] D.R. Sahu, Luoyi Shi, Ngai-Ching Wong and Jen-Chih Yao, Perturbed
iterative methods for a general family of operators: convergence theory and
applications, \textit{Optimization}, (2020).
https://doi.org/10.1080/02331934.\newline2020.1745798. Full access at:
https://www.tandfonline.com/.\medskip

[112] Esther Bonacker, Aviv Gibali and Karl-Heinz K\"{u}fer, Nesterov
perturbations and projection methods applied to IMRT, \textit{Journal of
Nonlinear and Variational Analysis}, Volume 4, Issue 1, Pages 63--86, 2020.
Available online at http://jnva.biemdas.com.
https://doi.org/10.23952/jnva.4.2020.1.06.\medskip

[113] A. Gibali, G.T. Herman and C. Schn\"{o}rr, Guest Editors,
Superiorization versus Constrained Optimization: Analysis and Applications,
Special Issue of J\textit{ournal of Applied and Numerical Optimization}
(JANO), Volume 2, Number 1, April 2020. Read the special issue on the
journal's website at: http:
//jano.biemdas.com/archives/category/volume-2-issue-1.\medskip

[114] A. Gibali, G.T. Herman and C. Schn\"{o}rr, Editorial: A special issue
focused on superiorization versus constrained optimization: analysis and
applications, \textit{Journal of Applied and Numerical Optimization} (JANO),
Vol. 2, pp. 1--2, 2020. Read it on the journal's website at:
http://jano.biemdas.com/issues/\newline JANO2020-1-1-1.pdf.\medskip

[115] Mokhtar Abbasi and Touraj Nikazad, Superiorization of block accelerated
cyclic subgradient methods, \textit{Journal of Applied and Numerical
Optimization} (JANO), Vol. 2, pp. 3--13, 2020. Read it on the journal's
website at: http://jano.biemdas.com/issues/JANO2020-1-2.pdf.\medskip

[116] Y. Censor, S. Petra and C. Schn\"{o}rr, Superiorization vs. accelerated
convex optimization: The superiorized / regularized least-squares case,
\textit{Journal of Applied and Numerical Optimization} (JANO), Vol. 2, pp.
15--62, 2020. Read it on the journal's website at:
http://jano.biemdas.com/issues/JANO2020-1-3.pdf.\medskip

[117] A. Gibali and M. Haltmeier, Superiorized regularization of inverse
problems, J\textit{ournal of Applied and Numerical Optimization} (JANO), Vol.
2, pp. 63--70, 2020. Read it on the journal's website at:
http://jano.biemdas.com/\newline issues/JANO2020-1-4.pdf.\medskip

[118] Gabor T. Herman, Problem structures in the theory and practice of
superiorization, \textit{Journal of Applied and Numerical Optimization}
(JANO), Vol. 2, pp. 71-76, 2020. Read it on the journal's website at: \newline
http://jano.biemdas.com/issues/JANO2020-1-5.pdf.\medskip

[119] T. Humphries, M. Loreto, B. Halter, W. O'Keeffe and L. Ramirez,
Comparison of regularized and superiorized methods for tomographic image
reconstruction, \textit{Journal of Applied and Numerical Optimization} (JANO),
Vol. 2, pp. 77--99, 2020. Read it on the journal's website at: 

http://jano.biemdas.com/issues/JANO2020-1-6.pdf.\medskip

[120] S. Reich and A.J. Zaslavski, Inexact orbits of set-valued nonexpansive
mappings with summable errors, \textit{Journal of Applied and Numerical
Optimization} (JANO), Vol. 2, pp. 101-107, 2020. Read it on the journal's
website at: http://jano.biemdas.com/issues/JANO2020-1-7.pdf.\medskip

[121] A.J. Zaslavski, Three extensions of Butnariu-Reich-Zaslavski theorem for
inexact infinite products of nonexpansive mappings, \textit{Journal of Applied
and Numerical Optimization} (JANO), Vol. 2, pp. 109--120, 2020. Read it on the
journal's website at: http://jano.biemdas.com/issues/JANO2020-1-8.pdf.\medskip

[122] Satoshi Takabe and Tadashi Wadayama, Deep unfolded multicast
beamforming. Preprint. (2020). https://arxiv.org/pdf/2004.09345.pdf. (This
work is related to [92] above).\medskip

[123] Elias S. Helou, Marcelo V. W. Zibetti and Gabor T. Herman, Fast Proximal
Gradient Methods for Nonsmooth Convex Optimization for Tomographic Image
Reconstruction, \textit{Sensing and Imaging}, accepted for publication.
(2020). https://arxiv.org/abs/2008.09720.\medskip

[124] Nikolai Janakiev, Superiorized Algorithms for Medical Image Processing:
Comparison of Shearlet-based Secondary Optimization Criteria, Marshall Plan
Scholarship Report, Salzburg University of Applied Sciences, Salzburg,
Austria, and The Graduate Center of the City University of New York, NY, USA,
September 2016, (63 pages). https://static1.squarespace.com/static/
559921a3e4b02c1d7480f8f4/t/5d4d476d5cae4e000125087e/1565345658569/\newline
Janakiev+Nikolai\_816.PDF.\medskip

[125] Lukas Kiefer, Stefania Petra, Martin Storath and Andreas Weinmann,
Multi-channel Potts-based reconstruction for multi-spectral computed
tomography, \textit{Inverse Problems}, 37(4):045004 (2021).\medskip

[126] Johan Alme, Gergely G\'{a}bor Barnaf\"{o}ldi, Rene Barthel, Vyacheslav
Borshchov, Tea Bodova, Anthony van den Brink, Stephan Brons, Mamdouh Chaar,
Viljar Eikeland, Grigory Feofilov, Georgi Genov, Silje Grimstad, Ola
Gr\o ttvik, H\aa vard Helstrup, Alf Herland, Annar Eivindplass Hilde, Sergey
Igolkin, Ralf Keidel, Chinorat Kobdaj, Naomi van der Kolk, Oleksandr
Listratenko, Qasim Waheed Malik, Shruti Mehendale, Ilker Meric, Simon Voigt
Nesb\o , Odd Harald Odland, G\'{a}bor Papp, Thomas Peitzmann, Helge Egil Seime
Pettersen, Pierluigi Piersimoni, Maksym Protsenko, Attiq Ur Rehman, Matthias
Richter, Dieter R\"{o}hrich, Andreas Tefre Samn\o y, Joao Seco, Lena
Setterdahl, Hesam Shafiee, \O istein Jelmert Skjolddal, Emilie Solheim, Arnon
Songmoolnak, \'{A}kos Sud\'{a}r, Jarle Rambo S\o lie, Ganesh Tambave, Ihor
Tymchuk, Kjetil Ullaland, H\aa kon Andreas Underdal, Monika
Varga-K\"{o}farag\'{o}, Lennart Volz, Boris Wagner, Fredrik Mekki Wider\o e,
RenZheng Xiao, Shiming Yang, Hiroki Yokoyama, A High-Granularity Digital
Tracking Calorimeter Optimized for Proton CT, \textit{Frontiers in Physics},
Vol. 8, article 568243, (20pp), (2020). DOI:10.3389/fphy.2020.568243.
https://www.frontiersin.org/\newline articles/10.3389/fphy.2020.568243/full
[Open access].\medskip

[127] Kaiwen Ma, Nikolaos V. Sahinidis, Sreekanth Rajagopalan, Satyajith
Amaran and Scott J. Bury, Decomposition in derivative-free optimization,
\textit{Journal of Global Optimization}, Vol. 81, 269--292 (2021).\newline
https://doi.org/10.1007/s10898-021-01051-w.\medskip

[128] Alexander J. Zaslavski, Nonsmooth Convex Optimization. In: Alexander J.
Zaslavski, \textit{The Projected Subgradient Algorithm in Convex
Optimization}. SpringerBriefs in Optimization. Springer, Cham., pp. 5--83,
2020. Available from the publisher at
https://doi.org/10.1007/978-3-030-60300-7\_2.\medskip

[129] Lukas Kiefer, \textit{Efficient Algorithms for Mumford-Shah and Potts
Problems}, Ph.D. Dissertation, Combined Faculty for the Natural Sciences and
Mathematics, Heidelberg University, Germany. (220 pp.) (2020).
DOI:https://\newline doi.org/10.11588/heidok.00029100. Available at:
heiDOK-The Heidelberg Document Repository.\medskip

[130] Nuttapol Pakkaranang, Poom Kumam, Yusuf I. Suleiman and Bashir Ali,
Bounded perturbation resilience of viscosity proximal algorithm for solving
split variational inclusion problems with applications to compressed sensing
and image recovery, \textit{Mathematical Methods in the Applied Sciences}, pp.
1--23, 2020. The journal's website is at:
https://onlinelibrary.wiley.com/doi/\newline10.1002/mma.7023.\medskip

[131] Peichao Duan and Xubang Zheng, Bounded perturbation resilience of
generalized viscosity iterative algorithms for split variational inclusion
problems, \textit{Applied Set-Valued Analysis and Optimization}, Vol. 2, pp.
49--61, 2020. DOI:10.23952/asvao.2.2020.1.04. Open Access on the journal's
website at: http://asvao.biemdas.com/archives/1051.\medskip

[132] Jochen Fink, Renato L.G. Cavalcante and Slawomir Stanczak, Multi-group
multicast beamforming by superiorized projections onto convex sets,
\textit{IEEE Transactions on Signal Processing,} Vol. 69, pp. 5708--5722,
2020. DOI:10.23952/asvao.2.2020.1.04.\medskip

[133] Sebastian Meyer, Marco Pinto, Katia Parodi and Chiara Gianoli, The
impact of path estimates in iterative ion CT reconstructions for clinical-like
cases, \textit{Physics in Medicine and Biology}, 66, 095007, 2021.\newline
https://doi.org/10.1088/1361-6560/abf1ff.\medskip

[134] Weizhe Han, Qianlong Wang and Weiwei Cai, Computed tomography imaging
spectrometry based on superiorization and guided image filtering,
\textit{Optics Letters}, Vol. 46, pp. 2208--2211, 2021.\newline
https://doi.org/10.1364/OL.418355.\medskip

[135] Howard Heaton, Samy Wu Fung, Aviv Gibali and Wotao Yin,
Feasibility-based Fixed Point Networks, 2021.\newline
https://arxiv.org/abs/2104.14090.\medskip

[136] Maria Guenter, \textit{Comparison of Algorithms for Solving the Least
Squares Problem with Applications in Computed Tomography}, Master of Science
Thesis, Faculty of Science, Department of Computer Science, Mathematics,
Physics and Statistics, The University of British Columbia (UBC), Okanagan,
BC, Canada. (94 pp) (2021).\newline
https://open.library.ubc.ca/cIRcle/collections/ubctheses/24/items/1.0398196.\medskip

[137] Chongyuan Shui, Yihong Wang, Weiwei Cai, and Bin Zhou, Linear
multispectral absorption tomography based on regularized iterative methods,
\textit{Optics Express}, Vol. 29, pp. 20889--20912, 2021, \newline
https://doi.org/10.1364/OE.421817.\medskip

[138] Kaiwen Ma, Nikolaos V. Sahinidis, Sreekanth Rajagopalan, Satyajith
Amaran and Scott J. Bury, Decomposition in derivative-free optimization,
\textit{Journal of Global Optimization}, Vol. 81, pp. 269--292, 2021.\newline
https://doi.org/10.1007/s10898-021-01051-w.\medskip

[139] Mark Brooke, \textit{Incorporation of Biological Factors in Radiation
Therapy Treatment Planning}, Ph.D. Dissertation, Department of Oncology,
Wolfson College, University of Oxford, UK. December 2020.\medskip

[140] Yair Censor, Keith E. Schubert and Reinhard W. Schulte, Developments in
mathematical algorithms and computational tools for proton CT and particle
therapy treatment planning, \textit{IEEE Transactions on Radiation and Plasma
Medical Sciences}, Vol. 6, pp. 313--324, (2022).\newline
DOI:10.1109/TRPMS.2021.3107322. Open Access on IEEEXplore.\medskip

[141] Yingying Li and Yaxuan Zhang, Bounded perturbation resilience of two
modified relaxed CQ algorithms for the multiple-sets split feasibility
problem, \textit{Axioms}, 10, 197, 2021 (22 pages).\newline
https://doi.org/10.3390/axioms10030197.\medskip

[142] Blake Edward Schultze, \textit{Essential Elements of Proton Computed
Tomography for Practical Applications}, Ph.D. Dissertation, Department of
Electrical and Computer Engineering, Baylor University, Waco, TX, USA. August
2021. Document Preview available at: ProQuest Dissertations
Publishing.\medskip

[143] Christina M. Sarosiek, \textit{Clinical Applications and Feasibility of
Proton CT and Proton Radiography}, Ph.D. Dissertation, Department of Physics,
Northern Illinois University, IL, USA. August 2021. Document Preview available
at: ProQuest Dissertations Publishing.\medskip

[144] T. Nikazad, M. Abbasi, L. Afzalipour and T. Elfving, A new step size
rule for the superiorization method and its application in computerized
tomography, \textit{Numerical Algorithms}, 2021.
https://doi.org/10.1007/s11075-021-01229-z.\medskip

[145] Jingyan Xu and Fr\'{e}d\'{e}ric Noo, Convex optimization algorithms in
medical image reconstruction -- in the age of AI, \textit{Physics in Medicine
and Biology}, 2021.
http://iopscience.iop.org/article/10.1088/1361-6560/ac3842.\medskip

[146] Maria Guenter, Steve Collins, Andy Ogilvy, Warren Hare and Andrew
Jirasek, Superiorization versus regularization: A comparison of algorithms for
solving image reconstruction problems with applications in computed
tomography, \textit{Medical Physics}, accepted for publication, 2021.\newline
https://doi.org/10.1111/mp.15373.\medskip

[147] Howard Wayne Heaton, \textit{Learning to Optimize with Guarantees},
Ph.D. Dissertation, Department of Mathematics, University of California, Los
Angeles, CA, USA. 2021. Dissertation available at: eScholarship
Publishing.\newline
https://escholarship.org/content/qt3274t029/qt3274t029.pdf.\medskip

[148] Ruiwen Xing, Thomas Humphries and Dong Si, Self-Attention Generative
Adversarial Network for Iterative Reconstruction of CT Images, 2022.
https://arxiv.org/abs/2112.12810.\medskip

[149] Qiao-Li Dong, Yeol Je Cho, Songnian He, Panos M. Pardalos and
Themistocles M. Rassias, The inertial Krasnosel'ski\u{\i}--Mann iteration, In:
\textit{The Krasnosel'ski\u{\i}--Mann Iterative Method: Recent Progress and
Applications}, SpringerBriefs in Optimization, Springer, Cham, Switzerland,
Chapter 5 (pp.59--73), 2022.
https://link.springer.com/chapter/10.1007/978-3-030-91654-1\_5.\medskip

[150] Jochen Fink, Renato L.G. Cavalcante and Slawomir Stanczak, Superiorized
Adaptive Projected Subgradient Method with Application to MIMO Detection,
2022. https://arxiv.org/abs/2203.01116.\medskip

[151] Yair Censor, Daniel Reem and Maroun Zaknoon, A generalized
block-iterative projection method for the common fixed point problem induced
by cutters, \textit{Journal of Global Optimization}, Vol. 84, pp. 967--987,
(2022). https://doi.org/10.1007/s10898-022-01175-7.\medskip

[152] Jochen Fink, Renato L.G. Cavalcante, Zoran Utkovski, and Slawomir
Stanczak, A Set-Theoretic Approach to Mimo Detection, in: \textit{IEEE
International Conference on Acoustics, Speech and Signal Processing
(ICASSP-2022)}, 2022, pp. 5328-5332.
https://ieeexplore.ieee.org/document/9746548.\medskip

[153] Idan Steinberg and Sanjiv Sam Gambhir (Inventors), Real-time
photoacoustic imaging using a precise forward model and fast iterative
inverse, Patent application number: 17/517401, Publication date: May 5, 2022.
See: https://portal.uspto.gov/pair/PublicPair. The full application appears at
the \textit{United States Patent and Trademark Office (USPTO)}, \newline
https://ppubs.uspto.gov/pubwebapp/static/pages/landing.html. This is related
to item [132] above.\medskip

[154] W. Cui, \textit{Research And Application Of Superiorization Algorithm
For Convex Optimization Problems}, Master Thesis, College of Science, Civil
Aviation University of China, Tianjin, P.R. China. Posted on: February 17,
2021. Thesis abstract and table of contents available at: Global Thesis.
https://www.globethesis.com/?t=2370330611468675.\medskip

[155] Oriol Garc\'{\i}a Llopis, Superiorizaci\'{o}n en Programaci\'{o}n Lineal
[English: Superiorization in Linear Programming], Bachelor Thesis, Department
of Mathematics, Faculty of Science, University of Alicante, Alicante, Spain.
Posted on: June 28, 2022. Thesis, in Spanish with English abstract, freely
available at: RUA. Repositorio Institucional de la Universidad de Alicante.
https://rua.ua.es/dspace/handle/10045/124555.\medskip

[156] Mahdi Mirzapour and Hossein Rabbani, Investigation on accelerated
ordered subsets image reconstruction techniques with superiorization
methodology, \textit{The European Physical Journal Plus}, Volume 137, Article
number: 791, 2022. https://doi.org/10.1140/epjp/s13360-022-02964-5.\medskip

[157] Francisco J. Arag\'{o}n-Artacho, Yair Censor, Aviv Gibali and David
Torregrosa-Bel\'{e}n, The superiorization method with restarted perturbations
for split minimization problems with an application to radiotherapy treatment
planning, \textit{Applied Mathematics and Computation}, Vol. 440, Article
127627 (2023). Open Access at:
https://doi.org/10.1016/j.amc.2022.127627.\medskip

[158] Xiangyu Zeng, Shuhao Xia, Kai Yang, Youlong Wu and Yuanming Shi,
Over-the-Air Computation for Vertical Federated Learning, in: \textit{2022
IEEE International Conference on Communications Workshops (ICC Workshops)},
pp. 788--793, 2022. https://ieeexplore.ieee.org/document/9814484.\medskip

[159] Florian Barkmann, Yair Censor and Niklas Wahl, Superiorization as a
novel strategy for linearly constrained inverse radiotherapy treatment
planning, Preprint, July 26, 2022. Available on arXiv at: \newline
https://arxiv.org/abs/2207.13187.\medskip

[160] Yiran Jia, Noah McMichael, Pedro Mokarzel, Brandon Thompson, Dong Si and
Thomas Humphries, Superiorization-inspired unrolled SART algorithm with U-Net
generated perturbations for sparse-view and limited-angle CT reconstruction,
\textit{Physics in Medicine \& Biology}, 67, 245004, (2022). Online at:
http://iopscience.iop.org/article/10.1088/1361-6560/aca513.\medskip

[161] Alexander J. Zaslavski, Superiorization with a projected subgradient
method, \textit{Journal of Applied and Numerical Optimization}, Vol. 4, pp.
291--298.(2022). Online at: http://jano.biemdas.com/archives/1357.\medskip

[162] Yair Censor, Superiorization: The asymmetric roles of
feasibility-seeking and objective function reduction, \textit{Applied
Set-Valued Analysis and Optimization}, accepted for publication, (2022).
\newline https://arxiv.org/abs/2212.14724.\medskip

[163] M\"{u}zeyyen Ert\"{u}rk and Ahmet Salkim, Superiorization and bounded
perturbation resilience of a gradient projection algorithm solving the convex
minimization problem, \textit{Optimization Letters}, (2023). \newline
https://link.springer.com/content/pdf/10.1007/s11590-022-01961-y.\medskip

[164] Yiji Wang, Cheng Zou, Dingzhu Wen and Yuanming Shi, Federated Learning
over LEO Satellite, \textit{2022 IEEE Globecom Workshops (GC Wkshps): Edge
Learning over 5G Mobile Networks and Beyond}, Rio de Janeiro, Brazil, (2022),
pp. 1652--1657.\newline
https://ieeexplore.ieee.org/abstract/document/10008719.\medskip

[165] Alexander J. Zaslavski, Two Convergence Results for Inexact Infinite
Products of Non-Expansive Mappings, \textit{Axioms}, 12, 88 (2023), (9pp).
Open access at: https://www.mdpi.com/2075-1680/12/1/88.\medskip

[166] Jochen Fink, \textit{Fixed Point Algorithms and Superiorization in
Communication Systems}, Ph.D. Dissertation, Fakult\"{a}t IV-Elektrotechnik und
Informatik, Technische Universit\"{a}t Berlin, Germany, 2022. Available at:
Technische Universit\"{a}t Berlin Publications. \newline
https://api-depositonce.tu-berlin.de/server/api/core/bitstreams/81740fd4-94f5-4c37-8804-b4cbb3dbf5ac/content.\medskip

[167] Adeolu Taiwo and Simeon Reich, Bounded perturbation resilience of a
regularized forward-reflected-backward splitting method for solving
variational inclusion problems with applications, \textit{Optimization},
(2023). \newline https://doi.org/10.1080/02331934.2023.2187664.\medskip

[168] Shan Yu, Yair Censor, Ming Jiang and Guojie Luo, Per-RMAP:
Feasibility-Seeking and Superiorization Methods for Floorplanning with I/O
Assignment, (2023). Presented at \textit{the International Symposium of EDA
(Electronics Design Automation) ISEDA-2023}, Nanjing, China, May 8-11, 2023.
Available on arXiv at: https://arxiv.org/abs/2304.06698.\medskip

[169] Jochen Fink, Renato L.G. Cavalcante, Zoran Utkovski, and Slawomir
Stanczak, Deep-unfolded adaptive projected subgradient method for MIMO
detection, in: \textit{IEEE International Conference on Acoustics, Speech and
Signal Processing (ICASSP-2023)}, 2023,\newline
DOI:10.1109/ICASSP49357.2023.10096706.\medskip

[170] Simeon Reich and Alexander J. Zaslavski, Convergence Results for Inexact
Iterates of Uniformly Locally Nonexpansive Mappings, \textit{Symmetry},
(2023), 15(5), 1084; https://doi.org/10.3390/sym15051084.\medskip

[171] Anitha Gopalan, O. Vignesh, R. Anusuya, K.P. Senthilkumar, V.S. Nishok,
T. Helan Vidhya and Florin Wilfred, Reconstructing the Photoacoustic Image
with High Quality using the Deep Neural Network Model, \textit{Contrast Media
\& Molecular Imaging}, Vol. 2023, Article ID 1172473, 13 pages, (2023).
https://doi.org/10.1155/2023/1172473.\medskip

[172] Samuel da Silva Oliveira, \textit{Parallelized Superiorization Method
for History Matching Problems Using Seismic Priors and Smoothness in Parts},
Ph.D. Dissertation, Department of Informatics and Applied Mathematics, The
Federal University of Rio Grande do Norte, Brazil, 2023. Dissertation
available at: https://repositorio.ufrn.br.\newline
https://repositorio.ufrn.br/bitstream/123456789/52424/1/\newline
Metodoparalelosuperiorizacao\_Oliveira\_2023.pdf.\medskip

[173] Navneet Agrawal, Renato L.G. Cavalcante, Masahiro Yukawa, and Slawomir
Sta\'{n}czak, Distributed Convex Optimization \textquotedblleft
Over-the-Air\textquotedblright\ in Dynamic Environments, 2023. Available on
arXiv at: \newline https://arxiv.org/abs/2307.04913.\medskip

[174] Yanni Guo and Weixia Wang, The scaled and generalized splitting
algorithms for the minimization problem of the sum of finite convex functions,
\textit{Journal of Nonlinear and Covnex Analysis}, Vol. 24, pp. 101--117,
2023. See the complete table of contents of the special issue at Yokohama
Publishers: http://www.yokohamapublishers.jp/online2/jncav24-1.html.\medskip

[175] Yan Tang and Zhihui Ji, Perturbation Resilience of Self-Adaptive
Step-Size Algorithms for Solving Split Variational Inclusion Problems and
their Applications, \textit{Numerical Functional Analysis and Optimization},
2023, DOI:10.1080/01630563.2023.2247615.\medskip

[176] Mokhtar Abbasi, Mahdi Ahmadinia and Ali Ahmadinia, A new step size
selection strategy for the superiorization methodology using subgradient
vectors and its application for solving convex constrained optimization
problems, \textit{IMA Journal of Numerical Analysis}, (2023), drad070.
https://doi.org/10.1093/imanum/drad070.\medskip

[177] Navneet Agrawal, Renato L.G. Cavalcante and Slawomir Sta\'{n}czak,
Distributed fixed-point algorithms for dynamic convex optimization over
decentralized and unbalanced wireless networks, 2024. Available on arXiv at:
https://arxiv.org/abs/2401.18030v1.\medskip

[178] Mokhtar Abbasi and Touraj Nikazad, Bounded perturbations resilient
iterative methods for linear systems and least squares problems:
operator-based approaches, analysis, and performance evaluation, \textit{BIT
Numerical Mathematics}, Vol. 64, Article number 15 (2024).
https://doi.org/10.1007/s10543-024-01015-y.\medskip

[179] Alexander J. Zaslavski, \textit{Solutions of Fixed Point Problems with
Computational Errors}, Springer Nature Switzerland AG, (2024).\medskip

[180] Alexander J. Zaslavski, Convergence of inexact orbits of nonexpansive
mappings in complete metric spaces, \textit{Communications in Optimization
Theory}, Vol. 2024 (2024), Article ID 4, pp. 1--10. \newline
https://cot.mathres.org/issues/COT20244.pdf.\medskip

[181] Steve Collins, Andy Ogilvy, Warren Hare, Michelle Hilts and Andrew
Jirasek, Iterative Image Reconstruction Algorithm Analysis for Optical CT
Radiochromic Gel Dosimetry, analysis, and performance evaluation,
\textit{Biomedical Physics \& Engineering Express}, accepted manuscript, open
access (2024).
https://iopscience.iop.org/article/10.1088/2057-1976/ad3afe.\medskip

[182] Rodrigo Hernang\'{o}mez, Jochen Fink, Renato L.G. Cavalcante, Zoran
Utkovski and Slawomir Sta\'{n}czak, Optimized Detection with Analog
Beamforming for Monostatic Integrated Sensing and Communication, 2024.
Available on arXiv at: https://arxiv.org/abs/2404.08455.\medskip

[***] F.J. Arag\'{o}n-Artacho, Y. Censor, A. Gibali and D.
Torregrosa-Bel\'{e}n, Corrigendum to: The superiorization method with
restarted perturbations for split minimization problems with an application to
radiotherapy treatment planning, \textit{Applied Mathematics and Computation},
Vol. 476, Article 128788 (2024). Available at:\newline
https://www.sciencedirect.com/science/article/pii/S0096300324002522.\medskip

[183] F.J. Arag\'{o}n-Artacho, W. Cai, Y. Censor, A. Gibali, C. Shui and D.
Torregrosa-Bel\'{e}n, Approaches to iterative algorithms for solving nonlinear
equations with an application in tomographic absorption spectroscopy,
\textit{Communications in Optimization Theory}, accepted for publication,
2024. Available on arXiv at: https://arxiv.org/abs/2405.08635.\medskip

[184] Shousheng Luo, Zhiting Liu, Yaofei Lu and Xue-Cheng Tai, Superiorized
iteration algorithm for CT image simultaneous reconstruction and segmentation,
\textit{Inverse Problems and Imaging}, DOI:10.3934/ipi.2024024, (2024).
Available at American Institute of Mathematical Sciences. \newline
https://www.aimsciences.org//article/doi/10.3934/ipi.2024024.\medskip

[185] S. Yu, Y. Censor and G. Luo, Floorplanning With I/O Assignment via
Feasibility-Seeking and Superiorization Methods, \textit{IEEE Transactions on
Computer-Aided Design of Integrated Circuits and Systems}, Vol. 44, pp.
317--330, (2024). Online at:\newline IEEEXplore.
https://ieeexplore.ieee.org/abstract/document/10543062.\medskip

[186] David Torregrosa Belen, \textit{Splitting Algorithms for Structured
Optimization: Theory and Applications}, Ph.D. Dissertation, Department of
Mathematics, University of Murcia, Spain, 2024. Dissertation available at:
Digitum: Repositorio Institucional de la Universidad de Murcia. \newline
https://digitum.um.es/digitum/handle/10201/142204.\medskip

[187] Touraj Nikazad, Perturbed conjugate gradient method: controlling
semi-convergence, enforcing non-negativity, and accelerating iterations in
tomographic imaging, \textit{Physica Scripta}, Vol. 99, 075034, (2024).
DOI:10.1088/1402-4896/ad5651. Available at: IOPscience. \newline
https://iopscience.iop.org/article/10.1088/1402-4896/ad5651/meta.\medskip

[188] Jan Schr\"{o}der, Yair Censor, Philipp S\"{u}ss and Karl-Heinz
K\"{u}fer, Immunity to increasing condition numbers of linear superiorization
versus linear programming, Preprint. Available on arXiv at: \newline
https://arxiv.org/abs/2407.18709.\medskip

[189] Chen JingRuo, Xu ShiJie, Liu HeCong, Huang JianQing, Liu YingZheng and
Cai WeiWei, Untrained neural network for linear tomographic absorption
spectroscopy, \textit{Science China Technological Sciences}, Vol. 67, (2024).
https://doi.org/10.1007/s11431-023-2629-2. Available on SpringerLink.\medskip

[190] Jon Henshaw, Aviv Gibali and Thomas Humphries, Plug-and-play
superiorization. Available on arXiv at:
https://arxiv.org/abs/2410.23401.\medskip

[191] Satoshi Shoji, Wataru Yata, Keita Kume and Isao Yamada, A
Discrete-Valued Signal Estimation by Nonconvex Enhancement of SOAV with cLiGME
Model, \textit{2024 Asia Pacific Signal and Information Processing Association
Annual Summit and Conference} (APSIPA ASC), GALAXY INTERNATIONAL CONVENTION
CENTER, Macau Cotai Macau, Macau (2024).
http://www.apsipa2024.org/welcome.html.\medskip

[192] K. Barshad, Y. Censor, W.M. Moursi, T. Weames and H. Wolkowicz, A
necessary condition for the guarantee of the superiorization method,
\textit{Optimization Letters}, https://doi.org/10.1007/s11590-025-02192-7
(2025). Open Access at Springer Nature Link.\medskip

[193] Satoshi Shoji, Wataru Yata, Keita Kume and Isao Yamada, An LiGME
Regularizer of Designated Isolated Minimizers -- An Application to
Discrete-Valued Signal Estimation, (2025). Available on arXiv at: \newline
https://arxiv.org/abs/2503.10126.\medskip

[194] Mingxia Zheng and Yanni Guo, Scaled forward-backward algorithm and the
modified superiorized version for solving the split monotone variational
inclusion problem, \textit{Optimization Eruditorum}, Vol. 1(1), pp. 56--74,
(2024). Freely available at: Tulipa Opera Scholarum.\medskip

[195] Alexander J. Zaslavski, Superiorization on solution sets of common fixed
point problems with countable families of maps, \textit{Pacific Journal of
Optimization}, Vol. 21, pp. 627--642, (2025).
https://doi.org/10.61208/pjo-2025-008.\medskip

[196] Damian Borys, Jan Gajewski, Tobias Becher, Yair Censor, Renata
Kope\'{c}, Marzena Rydygier, Angelo Schiavi, Tomasz Sk\'{o}ra, Anna Spaleniak,
Niklas Wahl, Agnieszka Wochnik and Antoni Ruci\'{n}ski, GPU-accelerated
FREDopt package for simultaneous dose and LETd proton radiotherapy plan
optimization via superiorization methods, \textit{Physics in Medicine and
Biology}, accepted for publication (2025). Available open access at: \newline
https://iopscience.iop.org/article/10.1088/1361-6560/ade841.\medskip

[197] Kay Barshad and Yair Censor, General Perturbation Resilient Dynamic
String-Averaging for Inconsistent Problems with Superiorization,
\textit{Journal of Optimization Theory and Applications}, Vol. 207, article
number 9, (2025). Available open access at:\newline
https://link.springer.com/article/10.1007/s10957-025-02763-9.\medskip

[198] Shan Yu, Haiyang Liu, Xinming Wei, Bizhao Shi and Guojie Luo, TACPlace:
Ultrafast Thermal-Aware Chiplet Placement with Feasibility Seeking,
\textit{GLSVLSI'25: Proceedings of the Great Lakes Symposium on VLSI 2025},
New Orleans, LA, USA. 30 June 2025-2 July 2025, pp. 600-605, (2025). Available
open access at: ACM Digital Library.\medskip

[199] Alexander J. Zaslavski, Superiorization technique with a projected
subgradient method for games, \textit{Pure and Applied Functional Analysis},
Vol. 8, Number 4, pp. 1211--1221, (2023).
http://yokohamapublishers.jp/online-p/Pafa/vol8/pafav8n4p1211.pdf.\medskip

[200] Qiao-Li Dong, Xiao-Huan Li and Songnian He, Outer perturbations of a
projection method and two approximation methods for the split equality
problem, \textit{Optimization}, Vol. 67, pp. 1429--1446, (2018). \newline
https://www.tandfonline.com/doi/full/10.1080/02331934.2018.1474470.\medskip

[201] Yuhui Nie, Mengyuan Wang, Yuheng Wang, Junjie Lin, Bingxin Liu, Tao Yin,
Zhipeng Liu and Shunqi Zhang, Superiorized Model-based Real-time Inversion for
Cross-sectional Magnetoacoustic Tomography Combined with Magnetic Induction,
\textit{Physics in Medicine and Biology}, accepted for publication (2025).
Available at:
https://iopscience.iop.org/article/10.1088/1361-6560/ae2cdd.\medskip

[202] Tobias Becher, Yair Censor and Niklas Wahl, LET based proton plan
calculation with superiorization and feasibility-seeking based dose mimicking,
\textit{International Journal of Particle Therapy}, Vol. 17, Supplement,
101015, (2025). Conference abstract SO061/\#149. Available open access at:
https://doi.org/10.1016/j.ijpt.2025.101015.\medskip

[203] Navneet Agrawal, Dynamic Distributed Optimization and Fixed-Point
Seeking via Over-the-Air Function Computation, Ph.D. Dissertation,
Fakult\"{a}t IV-Elektrotechnik und Informatik der Technischen Universit\"{a}t
Berlin, Berlin, Germany, 2024. Dissertation available at:
https://depositonce.tu-berlin.de/items/9f758380-2a12-43ac-b65a-bea8ec261a1d.\medskip

[204] Touraj Nikazad and Mokhtar Abbasi, Biased superiorization of steepest
descent: redefining the reconstruction target in noisy inverse problems,
\textit{Electronic Transactions on Numerical Analysis}, Vol. 65, pp. 93--109,
(2026).
https://etna.ricam.oeaw.ac.at/vol.65.2026/pp93-109.dir/pp93-109.pdf.\medskip

[205] Kay Barshad and Yair Censor, Strong convergence, perturbation resilience
and superiorization of Generalized Modular String-Averaging with infinitely
many input operators, \textit{Numerical Algorithms}, accepted for publication,
(2026).\medskip%

\[
To\text{ \ }be\text{ \ }continued...
\]

\end{document}